УДК 511.215

**БАЙБЕКОВ С.Н., ДУРМАГАМБЕТОВ А.А.**
г.Астана, Казахстан

## БЕСКОНЕЧНОСТЬ КОЛИЧЕСТВА ПРОСТЫХ ЧИСЕЛ-БЛИЗНЕЦОВ

Одной из главных проблем простых чисел, которая вот уже в течение более 2000 лет не имеет решения, является обозначенная 5-м Международным математическим конгрессом проблема о доказательстве бесконечности количества простых чисел-близнецов. В настоящее время эта проблема именуется как вторая проблема Ландау.

В 2013 году американский математик Чжан Итан (Yitang Zhang) из университета Нью-Гэмпшира доказал, что существует бесконечно большое количество пар простых чисел, расстояние между которыми не превышает 70 миллионов. Позже, Джеймс Мэйнард (James Maynard) улучшил результат до 600. В 2014 году проект Polymath под руководством Теренса Тао несколько улучшили последний метод, получив оценку в 246 [1].

Для решения этой проблемы нами был предложен новый метод, который эмпирическим способом позволяет оценить бесконечность простых чисел-близнецов [2]. В настоящей работе предлагается другое простое доказательство, которое дает более корректный результат о бесконечности близнецов.

Сначала для удобства введем следующее обозначение. Как известно последовательное умножение натуральных чисел называется факториалом: $\prod_{i=1}^{n} i = n!$ В дальнейшем последовательное умножение простых чисел будет встречаться довольно часто, поэтому для таких случаев используем следующее обозначение:

$$2*3*5*7*11*...*p_n = \prod_{i=1}^{n} p_i = p_n!'$$

Здесь $p_i$ - простое число с порядковым номером $i$. Комбинация знаков $p_n!'$ означает последовательное произведение исключительно простых чисел от 2 до $p_n$. Назовем его **специальным факториалом простого числа $p_n$**. Например, $p_4!'$ - специальный факториал простого числа $p_4 = 7$ или $p_4!' = 7!' = 2*3*5*7 = 210$.

**Матрицы простых чисел.**

Представим множество натуральных чисел в виде семейства матриц $A_k$ с элементами $a(k,i,j)$, где $i$-порядковый номер строк, а $j$-порядковый номер столбцов матрицы, $k$-порядковый номер матрицы $A_k$.

При этом максимальное количество строк матрицы $A_k$ должно быть равно специальному факториалу $p_k!'$, т.е. $i_{k,max} = p_k!'$. Это означает, что каждой матрице с порядковым номером $k$ должен соответствовать определенный набор последовательности простых чисел: $p_1, p_2, p_3, \ldots, p_k$ (заметим, последним простым числом, которое соответствует этой матрице, является $p_k$). Количество столбцов может быть сколь угодно большим вплоть до бесконечности.

Здесь и в дальнейшем полагаем, что нам не известно никакое простое число. Простые числа будут генерироваться по ходу построения матриц $A_k$.

Сначала покажем, как формируется матрица $A_1$. Для этого рассмотрим ряд натуральных чисел от 1 до бесконечности *(рис.1, $A_0$)*. В этом ряду после числа 1 следует число 2. Стало быть, число 2 делится только на 1 и на 2. Следовательно, первым простым числом является 2, т.е. $p_1 = 2$. Тогда первая матрица $A_1$, построенная с учетом первого простого числа, имеет всего 2 строки ($i_{k,max} = i_{1,max} = p_1!' = 2!' = 2$). Количество столбцов бесконечно *(рис.1, $A_1$-а)*. В первой строке матрицы $A_1$ расположены числа 2, 4, 6, ... Эти числа образуют арифметическую прогрессию. Первый член, а также разность (шаг) этой прогрессии равны 2, т.е. члены первой строки соответственно равны

$$a(k,i,j) = a(1,1,j) = 2 + 2*(j-1),\text{ где } j=1, 2, 3,\ldots$$

Числа, расположенные во второй строке новой матрицы, тоже образуют арифметическую прогрессию. Первый член и разность этой прогрессии соответственно равны 3 и 2, т.е.

$$a(k,i,j) = a(1,2,j) = 3 + 2*(j-1),\text{ где } j=1, 2, 3,\ldots$$

Как уже определили, число 2 является простым числом. Поэтому все числа, кратные 2, являются составными числами. В силу этого все числа, за исключением 2, расположенные в первой строке рассматриваемой матрицы *(рис.1, $A_1$-в)*, для наглядности окрашиваются темным цветом. Тем самым при помощи матрицы $A_1$ определяются все составные числа, которые должны делиться на 2.

Отсюда следует, что число 1 не является простым числом, иначе все числа, кратные единице, были бы составными числами. Число 1 также не является составным числом, так как оно не делится на другие числа. Именно поэтому число 1 в этой и других матрицах располагается отдельно в верхнем левом уголке.

**Алгоритм преобразования матрицы с одного вида в другой вид**.

Из матрицы $A_1$ *(рис.1, $A_1$-в)* видно, что после числа 2 не окрашенным числом является число 3, оно не делится на 2 и в силу этого является вторым простым числом, т.е. $p_2 = 3$.

Поэтому матрицу $A_1$ *(рис. 1, $A_1$-в)* преобразуем в следующую матрицу $A_2$ *(рис.1, $A_2$-а)*. Максимальное количество строк этой матрицы должно быть равно специальному факториалу второго простого числа $p_2 = 3$, т.е. $i_{k,max} = i_{2,max} = p_2!' = 3!' = 6$. Количество столбцов, как в первом случае, может быть сколь угодным.

Для преобразования матрицы с одного вида в другой вид используется несложный способ, реализация которого заключается в простом переносе чисел из определенных строк и столбцов исходной матрицы в соответствующие строки и столбцы новой матрицы. Например, для формирования первого столбца матрицы $A_2$ сначала числа 2 и 3, расположенные в первом столбце исходной матрицы $A_1$, перенесутся в первую и вторую строку новой матрицы $A_2$. Затем числа 4 и 5 из первой матрицы перенесутся в 3-ю и 4-ю строку новой матрицы. Затем числа 6 и 7 также перенесутся в 5-ю и 6-ю строку новой матрицы. На этом завершается формирование первого столбца матрицы $A_2$.

Для формирования второго столбцы матрицы $A_2$ аналогичным образом числа (8 и 9), (10 и 11) и (12 и 13) попарно перенесутся в 2-й столбец матрицы $A_2$. Далее аналогичным образом создаются и другие столбцы.

Заметим, что в новой матрице (*рис.1, $A_2$-а*) все числа, расположенные в 3-й и 5-й строках, окрашены темным цветом, так как они благодаря матрице $A_1$ уже были определены как составные числа.

В новой матрице все числа, расположенные в каждой строке, также как в случае матрицы $A_1$, образуют арифметическую прогрессию, которая в общем случае имеет следующий вид:

$$a(k, i, j) = (i + 1) + p_k!'(j - 1), \text{ где } j = 1, 2, \ldots, \infty; \quad i = 1, 2, \ldots, p_k!' \quad (1)$$

Разностью этой арифметической прогрессии является $p_k!'$. Выражение (1) для матрицы $A_2$, в частности для ее 2-й строки, имеет следующий вид:

$$a(2,2,j) = 3 + 6(j - 1), \text{ где } j = 1, 2, \ldots, \infty; \quad i = 1, 2, \ldots, p_2!'$$

Как видно, все числа этой строки делятся на 3, следовательно, они (за исключением числа 3) являются составными числами. В силу этого они в матрице $A_2$ *(рис.1, $A_2$-в)* окрашены темным цветом. Что касается чисел, расположенные в 5-й строке, они тоже делятся на 3. Но, эти числа, как было сказано выше, при рассмотрении матрицы $A_1$ уже были определены как составные числа. Множество чисел 4-й строки (также 6-й строки) тоже образует арифметическую прогрессию. Но среди них имеются и простые и составные числа. Поэтому числа этих строк пока не будут окрашены. Заметим, что в темный цвет окрашиваются те строки, в которых находятся только составные числа.

Обратим внимание, что здесь и во всех последующих рисунках буквой – *а* будут обозначены те матрицы (например, $A_1$-а, $A_2$-а, $A_3$-а), которые получаются после преобразования предыдущей матрицы. Буквой – *в* будут обозначены те матрицы (например, $A_1$-в, $A_2$-в, $A_1$-в), которые получаются после переработки уже преобразованных матриц.

**Ряд натуральных чисел**

| 1 | 2 | 3 | 4 | 5 | 6 | 7 | 8 | 9 | 10 | ... | ... | ... | ... |
|---|---|---|---|---|---|---|---|---|----|-----|-----|-----|-----|

A0

**A1-а**

| | | Номера столбцов | | | |
|---|---|---|---|---|---|
| | 1 | 1 | 2 | 3 | 4 | ... |
| | 1 | | 2 | 4 | 6 | 8 | ... |
| | 2 | | 3 | 5 | 7 | 9 | ... |

**A1-в**

| | | Номера столбцов | | | |
|---|---|---|---|---|---|
| | 1 | 1 | 2 | 3 | 4 | ... |
| | 1 | | 2 | 4 | 6 | 8 | ... |
| | 2 | | 3 | 5 | 7 | 9 | ... |

**A2-а**

| | | | Номера столбцов | | | |
|---|---|---|---|---|---|---|
| | | 1 | 1 | 2 | 3 | 4 | ... |
| с | 1 | | 2 | 8 | 14 | 20 | ... |
| т | 2 | | 3 | 9 | 15 | 21 | ... |
| р | 3 | | 4 | 10 | 16 | 22 | ... |
| о | 4 | | 5 | 11 | 17 | 23 | ... |
| к | 5 | | 6 | 12 | 18 | 24 | ... |
| и | 6 | | 7 | 13 | 19 | 25 | ... |

**A2-в** (аналогичная таблица)

**A3-а**

| | | | Номера столбцов | | | |
|---|---|---|---|---|---|---|
| | | 1 | 1 | 2 | 3 | 4 | ... |
| | 1 | | 2 | 32 | 62 | 92 | ... |
| н | 2 | | 3 | 33 | 63 | 93 | ... |
| о | 3 | | 4 | 34 | 64 | 94 | ... |
| м | 4 | | 5 | 35 | 65 | 95 | ... |
| е | 5 | | 6 | 36 | 66 | 96 | ... |
| р | 6 | | 7 | 37 | 67 | 97 | ... |
| а | 7 | | 8 | 38 | 68 | 98 | ... |
| | 8 | | 9 | 39 | 69 | 99 | ... |
| | 9 | | 10 | 40 | 70 | 100 | ... |
| с | 10 | | 11 | 41 | 71 | 101 | ... |
| т | 11 | | 12 | 42 | 72 | 102 | ... |
| р | 12 | | 13 | 43 | 73 | 103 | ... |
| о | 13 | | 14 | 44 | 74 | 104 | ... |
| к | 14 | | 15 | 45 | 75 | 105 | ... |
| | 15 | | 16 | 46 | 76 | 106 | ... |
| | 16 | | 17 | 47 | 77 | 107 | ... |
| | 17 | | 18 | 48 | 78 | 108 | ... |
| | 18 | | 19 | 49 | 79 | 109 | ... |
| | 19 | | 20 | 50 | 80 | 110 | ... |
| | 20 | | 21 | 51 | 81 | 111 | ... |
| | 21 | | 22 | 52 | 82 | 112 | ... |
| | 22 | | 23 | 53 | 83 | 113 | ... |
| | 23 | | 24 | 54 | 84 | 114 | ... |
| | 24 | | 25 | 55 | 85 | 115 | ... |
| | 25 | | 26 | 56 | 86 | 116 | ... |
| | 26 | | 27 | 57 | 87 | 117 | ... |
| | 27 | | 28 | 58 | 88 | 118 | ... |
| | 28 | | 29 | 59 | 89 | 119 | ... |
| | 29 | | 30 | 60 | 90 | 120 | ... |
| | 30 | | 31 | 61 | 91 | 121 | ... |

**A3-в** (аналогичная таблица)

*Рис.1.* **Матрицы простых чисел**

Следует сказать, что в ходе обработки 2-й матрицы все числа, которые должны делиться на 3, выявляются окончательно, и они тоже соответствующим образом перекрашиваются темным цветом.

Из матрицы $A_2$ *(рис.1, $A_2$-в)* видно, что после чисел 2 и 3 неокрашенным числом является число 5. Следовательно, оно является третьим простым числом $p_3 = 5$. Поэтому 2-я матрица *(рис.1, $A_2$-в)* преобразуется в 3-ю матрицу $A_3$ *(рис.1, $A_3$-а)*. Для этого используется процедура, аналогичная процедуре, которая была использована для преобразования матрицы $A_1$ в матрицу $A_2$.

Например, для формирования первого столбца матрицы $A_3$ сначала числа (2÷7), расположенные в первом столбце исходной матрицы $A_2$, перенесутся в (1-ю ÷ 6-ю) строку новой матрицы $A_3$. Затем числа (8÷13) из второго столбца матрицы $A_2$ перенесутся в (7-ю ÷ 12-ю) строку новой матрицы. После этого числа (14 ÷19) перенесутся в (13-ю ÷ 18-ю) строку новой матрицы. Затем числа (20÷25) перенесутся в (19-ю ÷ 24-ю) строку новой матрицы и наконец, числа (26÷31) перенесутся в (25-ю ÷ 30-ю) строку новой матрицы. На этом завершается формирование первого столбца матрицы $A_3$.

Для формирования второго столбцы матрицы $A_3$ аналогичным образом числа (32÷37), (38÷43), (44÷49), (50÷55) и (56÷61), расположенные в матрице $A_2$, поэтапно перенесутся в 2-й столбец матрицы $A_3$. После этого аналогичным образом создаются и другие столбцы.

Максимальное количество строк 3-й матрицы должно быть равно специальному факториалу третьего простого числа $p_3 = 3$, то есть
$$i_{k,max} = i_{3,max} = p_3!' = 5!' = 30.$$
Для матрицы $A_3$ выражение (1) имеет следующий вид:

$$a(3, i, j) = (i + 1) + 30(j - 1), \text{где } j = 1, 2, \ldots, \infty; \ i = 1, 2, \ldots, p_3!'$$

Из этого выражения получим, что все числа, расположенные в 4-й строке делятся на 5, а те числа, которые расположены в 24-й строке, также делятся на число 5. В связи с этим все они соответствующим образом переходят в ряды составных чисел и перекашиваются в темный цвет (за исключением простого числа 5). Здесь в случае 3-й матрицы *(рис.1, $A_3$-в)* также следует отметить, что все числа, которые должны делиться на 5, выявляются окончательно и перекрашиваются в темный цвет (например, в частности строка с порядковым номером 24).

Обратим внимание, что во всех случаях расположение окрашенных и неокрашенных чисел в пределах одного столбца любой рассматриваемой матрицы не поддается никакой строгой закономерности. Но картина взаимного расположения этих чисел в пределах одного столбца повторяется с идеальной точностью в последующих столбцах (начиная со 2-го столбца). Такая закономерность повторения картин по столбцам получается тогда, когда каждая предыдущая матрица $A_n$ с количеством строк равным $p_n!'$

преобразуется в следующую матрицу $A_{n+1}$ с количеством строк равным $p_{n+1}!'$.

В матрице $A_3$ *(рис.1, $A_3$-в)* после чисел 2, 3, и 5 не окрашенным числом является число 7. Следовательно, оно является четвертым простым числом $p_4 = 7$. Теперь зная четвертое простое число 7, матрицу $A_3$ аналогичным образом можно преобразовать в следующую 4-ю матрицу $A_4$. У этой матрицы максимальное количество строк должно быть равно специальному факториалу простого числа 7, т.е. 7!'=210.

В этом случае проводя ряд аналогичных операций, как в предыдущих случаях, окончательно можно выявить множество всех составных чисел, которые должны делиться на число 7. Аналогично можно построить и другие матрицы.

Теперь при помощи матриц $A_k$ попробуем определить количество простых чисел-близнецов.

**Бесконечность количества простых чисел-близнецов.**

Сначала введем ряд определений:

*Определение 1*. Если в какой-нибудь строке матрицы имеются только составные числа, то для наглядности она окрашивается в темный цвет и для удобства называем ее <u>окрашенной строкой.</u>

*Определение 2*. Если в какой-нибудь строке имеются и простые и составные числа, то она не окрашивается, такие строки для удобства называются <u>неокрашенными строками.</u>

*Определение 3*. Если первое число строки является неокрашенным, а остальные числа являются окрашенными, то это неокрашенное число является простым числом, а остальные – составными числами.

*Определение 4*. Если разность порядковых номеров двух соседних неокрашенных строк равна 2, то таких строк называем <u>парой строк-близнецов</u> или <u>строками-близнецами.</u> Для чисел, расположенных в разных строках, но в одном столбце пар строк-близнецов, всегда выполняется равенство: $|a(k,i,j) - a(k, i \mp 2, j)| = 2$

*Определение 5*. Если порядковый номер какой-нибудь неокрашенной строки различается от порядкового номера ближайшей неокрашенной строки больше чем на 2, то такая строка называется <u>одинокой строкой.</u>

Из этих определений следуют, что простые числа-близнецы могут находиться <u>только</u> в строках-близнецах.

Цель настоящей работы заключается в определении общего количества простых чисел-близнецов. Поэтому в дальнейшем основное внимание будет <u>уделено</u> парам строк-близнецов.

**Теорема 1.** *Количество пар строк-близнецов, находящихся в матрице $A_k$, монотонно растет с ростом порядкового номера k матрицы, причем в каждой строке любой пары строк-близнецов находятся бесконечное количество простых чисел.*

Как известно, все простые числа-близнецы могут находиться только в парных строках-близнецах. При этом если в какой-то момент, например, при рассмотрении матрицы $A_k$, исчезнут все пары строк-близнецов, то, очевидно, в последующих матрицах их вообще не будут. В таком случае это означает, что количество простых чисел-близнецов должно быть ограниченным.

Проанализируем – возможен ли такой случай и за одно докажем теорему 1.

**Доказательство теоремы 1.**

Предположим, что в какой-то матрице имеется всего одна единственная пара строк-близнецов (например, как в случае *рис.1*, $A_2$-*в*). При этом имеется основание полагать, что в ходе дальнейшего преобразования данной матрицы в последующие матрицы, рассматриваемые пары строк-близнецов могут исчезнуть. Но, в самом деле, происходить наоборот. При преобразовании этой матрицы в следующую матрицу количество пар строк-близнецов, как было показано выше, становится **<u>больше</u>**.

Например, в матрице $A_2$ имеется всего одна единственная пара строк-близнецов *(рис.1, $A_2$-в)*. Из (1) следует, что члены арифметической прогрессии, которые находятся в строках этой единственной пары строк-близнецов, определяются выражением:

$$(6 \mp 1) + 6(j-1) = (p_2!' \mp 1) + p_2!'(j-1), \text{ где } j = 1, 2, \ldots, \infty \qquad (2)$$

Здесь знак «-» соответствует верхней, а знак «+» - нижней строке пары строк-близнецов.

Но эта единственная пара строк при преобразовании данной матрицы в матрицу $A_3$ порождает 5 ($p_3 = 5$) новых пар строк. То есть исходная единственная пара строк-близнецов разгруппируется по 5 новым парам строк.

Множество чисел, находящихся в каждой строке новых 5 пар строк матрицы $A_3$, тоже образует арифметическую прогрессию с постоянной $p_3!' = 30$ и определяется выражением:

$$(6m \mp 1) + 30(j-1) = (p_2!'m \mp 1) + p_3!'(j-1), \qquad (3)$$
где $j = 1, 2, 3, \ldots, \infty; \qquad m = 1, 2, \ldots, p_3.$

Из (3) получим, что если при каком-то значении $m = 1, 2, \ldots, p_3$ выполняется равенство

$$\frac{p_2!'m \mp 1}{p_3} = \text{целое число}, \qquad (4)$$

то все числа данной строки делятся на $p_3$ без остатков, следовательно, они являются составными числами. В самом деле, известно, что в пределах интервала $0 < m < p_3$ равенство (4) относительно параметра $m$ имеет единственное решение [3], [4], [5]. Например, равенство (4) для случая $(p_2!'m - 1)$ выполняется при $m = 1$, а для случая $(p_2!'m + 1)$ – при $m = 4$. Т.е. при $m = 1$ и $m = 4$ рассматриваемая пара строк не является парой строк-близнецов, и соответствующая строка, для которой выполняется равенство (4), перекрашивается в темный цвет. В результате из 5 вновь образованных пар строк только 3 пары являются парами строк-близнецов.

Все числа каждой строки вновь образованных 3 пар строк-близнецов, как было сказано выше, образуют арифметическую прогрессию, причем в каждой из них первый член и разность прогрессии являются взаимно простыми, т.е.:
$$(p_2!'m \mp 1, p_3!') \equiv 1, \text{где } m = 2, 3 \text{ и } 5, m \neq 1, m \neq 4$$

В силу этого из теоремы Дирихле для простых чисел в арифметической прогрессии следует, что в каждой строке этих 3 пар строк-близнецов имеется бесконечное количество простых чисел.

Теперь рассмотрим преобразование матрицы $A_3$ в матрицу $A_4$. В этом случае каждая пара строк-близнецов матрицы $A_3$ порождает по 7 ($p_4 = 7$) новых пар строк, итого в новой матрице $A_4$ образуются 21 новых пар строк. Значения чисел, находящихся в строках этих пар, определяются выражением:
$$[p_2!'p_i + p_3!'(m-1) \mp 1] + p_4!'(j-1), \qquad (5)$$

где $j = 1, 2, 3, \ldots, \infty$; $m = 1, 2, \ldots, p_4$; $i = 1, 2, 3$.

Из (5) видно, что множество чисел, находящихся в каждой строке вновь образованных 21 пар строк, в отдельности образует арифметическую прогрессию с разностью $p_4!'$ и с первыми членами $[p_2!'p_i + p_3!'(m-1) \mp 1]$.

Теперь рассмотрим, какие из этих 21 пар строк матрицы $A_4$ являются парами строк-близнецов. Для этого анализируем делимость первых членов вышесказанных арифметических прогрессий на $p_4 = 7$. При этом для удобства и наглядности рассмотрим случай, когда $i = 3$, но одновременно подразумеваем случаи, когда $i = 1$ и $i = 2$. Тогда из (5) получим, что значения чисел, находящихся в 7 новых парах строк матрицы $A_4$, порождаемые последней парой строк-близнецов матрицы $A_3$, определяются выражением:
$$p_3!'m \mp 1 + p_4!'(j-1),$$
где $j = 1, 2, 3, \ldots, \infty$; $m = 1, 2, \ldots, p_4$.

Рассмотрим делимость первых членов $p_3!'m \mp 1$ рассматриваемой арифметической прогрессии на значение $p_4$, т.е. выполняемость равенства:
$$\frac{p_3!'m \mp 1}{p_4} = \text{целое число}$$

В этом случае так же как в случае (4), получим, что в пределах интервала $0 < m < p_4$ данное равенство относительно параметра $m$ (для случая $p_3!'m - 1$, а также и для случая $p_3!'m + 1$) имеет единственное решение. Следовательно, в данном случае 2 пары рассматриваемых строк перестают быть парой строк-близнецов. Если к этому дополнительно рассмотрим случаи, когда $i = 1$ и $i = 2$, то получим, что в итоге 6 пар строк из 21 вновь образованных пар строк перестают быть парами строк-близнецов, и соответствующие строки, как было показано выше, перекрашиваются в темный цвет. В результате количество новых пар строк-близнецов в матрице $A_4$ будет равно 15.

Причем, все первые члены и разность арифметической прогрессии, образованных из чисел, находящихся в каждой строке вновь образованных 15 пар строк-близнецов матрицы $A_4$, являются взаимно простыми, т.е.:
$$([p_2!'p_i + p_3!'(m-1) \mp 1], p_4!') \equiv 1, \qquad (6)$$
где $i = 1, 2, 3;\ m \in (1, p_4)$ и $m \neq 5, 6, 9, 10, 17, 18$.

В силу этого из теоремы Дирихле для простых чисел в арифметической прогрессии следует, что в каждой строке вновь образованных 15 пар строк-близнецов имеются бесконечное количество простых чисел.

Если рассмотрим дальнейшие аналогичные преобразования матриц в последующие матрицы, например, матрицы $A_{k-1}$ в матрицу $A_k$, то каждый раз убеждаемся, что любая пара строк-близнецов исходной матрицы в новой матрице порождает $p_k$ штук новых пар строк. Причем 2 пары из них не будут парами строк-близнецов, и соответствующие строки переходят в ряды окрашенных и одиночных строк. Из этого получим, что в любой матрице $A_k$ общее количество строк ($i_{k,max}$) и общее количество пар строк-близнецов ($m_k$) соответственно равны:
$$i_{k,max} = i_{k-1,max} * p_k \quad \text{и} \quad m_k = m_{k-1}(p_k - 2) \qquad (7-1)$$
или
$$i_{k,max} = p_k!' \quad \text{и} \quad m_k = (p_k - 2)!', \quad k \geq 2 \qquad (7-2)$$

где $k$ – порядковый номер матрицы $A_k$ и/или простого числа $p_k$ и
$$(p_k - 2)!' = (p_2 - 2)(p_3 - 2) * \ldots * (p_k - 2) =$$
$$= 1 * 3 * 5 * 9 * 11 * \ldots * (p_k - 2)$$

Из (7) следует, что количество пар строк-близнецов по мере перехода в последующие матрицы, т.е. с ростом порядкового номера $k$ матрицы $A_k$ монотонно растет. С другой стороны, множество чисел, находящихся в каждой строке этих $m_k$ пар строк-близнецов, образуют арифметическую прогрессию. Причем первый член и разность каждой из этих прогрессии являются взаимно простыми. В силу этого из теоремы Дирихле следует, что в каждой строке любой пары строк-близнецов имеется бесконечное количество простых чисел.

**Теорема 1 доказана.**

Как было показано в (7), количество пар строк-близнецов по мере перехода в последующие матрицы, будет непрерывно расти. Но, при этом, количество общих строк каждой последующей матрицы растет как специальный факториал $p_k!'$. В силу этого **плотность** пар строк-близнецов по матрицам непрерывно уменьшается, так как с ростом $k$ монотонно уменьшается отношение $\frac{(p_k-2)!'}{p_k!'}$.

**Теорема 2.** *В любой паре строк-близнецов любой матрицы $A_k$ имеются простые числа-близнецы.*

Как было показано выше, все простые числа-близнецы могут находиться **только** в строках-близнецах. Но, могут ли быть случаи, когда в какой-то паре строк-близнецов все простые числа без исключения располагаются перекос на перекос, т.е. в силу этого никакая пара из двух простых чисел не могут находиться в одном столбце. Тогда в этой паре строк-близнецов не будет ни одной пары простых чисел-близнецов. Если во всех парах строк-близнецов данной матрицы случится такое, то в этой и во всех последующих матрицах больше не будут простых чисел-близнецов. Тогда однозначно можно сказать, что количество простых чисел-близнецов должно быть ограниченным.

Проанализируем этот случай и заодно докажем теорему 2.

**Доказательство теоремы 2.**
Рассмотрим матрицу $A_2$, в которой имеется всего одна единственная пара строк-близнецов *(рис.1, $A_2$-в)*. С другой стороны, как было показано выше, простые числа-близнецы могут находиться **только** в парах строк-близнецов. Это означает, что все существующие простые числа-близнецы находятся в этой единственной паре строк-близнецов.

Простой анализ показывает, что пары простых чисел-близнецов подчиняются некоторым простым правилам. В частности, последняя цифра любого простого числа (кроме 2 и 5) не может быть четным числом, она также не может быть равна 5. Это означает, что последними цифрами первого числа и второго числа любой пары близнецов соответственно должны быть (1 и 3) и (7 и 9), а также (9 и 1). Следовательно, множество простых чисел-близнецов по этим признакам должны разделиться на 3 подмножество. В самом деле, когда формируем матрицу $A_3$, единственная пара строк-близнецов матрицы $A_2$ порождает в матрице $A_3$ новые три пары строк-близнецов. Причем последняя цифра каждого числа, находящихся в строках отдельно выбранной пары строк-близнецов, соответственно равны (1 и 3) и (7 и 9), а также (9 и 1), *это наглядно показано в рис.1, $A_3$-в*.

Следовательно, то, что в каждой паре строк-близнецов матрицы $A_3$ имеется множество простых чисел-близнецов, причем сортированное с учетом значения последней цифры, не вызывает сомнения. То есть в парах строк-близнецов матрицы $A_3$ всеобщего перекоса простых чисел не будет.

Теперь рассмотрим следующую матрицу $A_4$. Проведенный нами анализ показывает, что все члены каждой пары строк-близнецов этой матрицы являются более упорядоченными, нежели в случае матрицы $A_3$. Например, последние 2 цифры каждого числа арифметической прогрессии, находящейся в любой строке любой пары строк-близнецов матрицы $A_4$, образуют циклически растущую последовательность типа: 21, 31, 41,…, 81, 91, 01, 11, 21, 31, …

При этом полагаем, что среднее расстояние _по столбцам_ между ячейками, где располагаются соседние простые числа _в одной_ паре строк-близнецов матрицы $A_4$, должно быть меньше среднего расстояния ячеек, где располагаются соседние простые числа _в одной_ паре строк-близнецов матрицы $A_3$.

Если это так, то в каждой паре строк-близнецов матрицы $A_4$ неизбежно существуют простые числа-близнецы. Так как в этом случае из-за тесноты вышеназванный перекос простых чисел если будет в матрице $A_4$, то будет в меньшей степени, чем в случае матрицы $A_3$. По крайней мере, в матрице $A_4$ всеобщего перекоса простых чисел не будет.

Чтобы убедиться в этом сначала вводим новый параметр $d_{i,i+1}$ - _расстояние по столбцам между ячейками_, где располагаются соседние простые числа с порядковыми номерами $i$ и $i+1$ _(рис.2)_:
$$d_{i,i+1} = |M_{i+1} - M_i|, \qquad (8)$$
где $M_i$ – порядковый номер столбца, в котором находится ячейка -го простого числа. При этом если два соседнего простого числа располагаются в двух примыкающих друг к другу ячейках по одной строке (т.е. по горизонтали, как это показано в _рис.2, а_), то расстояние между этими ячейками равно 1:
$$d_{i,i+1} = |M_{i+1} - M_i| = 1$$
С другой стороны, если два соседнего простого числа расположены в двух соседних и примыкающих друг к другу ячейках, находящихся в одном столбце (т.е. по вертикали, как это показано в _рис.2, в_), то расстояние между этими ячейками равно нулю:
$$d_{i,i+1} = |M_i - M_i| = 0$$
В этом случае эти два простого числа являются близнецами.

Теперь для примера рассмотрим фрагмент одной пары строк-близнецов произвольной матрицы. Пусть в этом фрагменте содержатся $n$ простых чисел, как это показано в _рис.2,с_. В этом рисунке светлыми кружочками обозначены ячейки, в которых располагаются составные числа, а темными кружочками обозначены ячейки простых чисел. Тогда вышеупомянутое _среднее расстояние по столбцам_ $d_\text{ср}$ между ячейками, где располагаются соседние простые числа, будет равно:

$$d_{\text{ср}} = \frac{\sum_{i=1}^{n-1} d_{i,i+1}}{n-1} \qquad (9)$$

В этом случае, как было показано выше, расстояние по столбцам между ячейками, где располагаются соседние простые числа с порядковыми номерами (1 и 3) и (4 и 6), равно 1, т.е.: $d_{1,3} = |M_3 - M_2| = 1$; $d_{4,6} = |M_6 - M_5| = 1$.

С другой стороны расстояние по столбцам между ячейками соседних простых чисел с порядковыми номерами (1 и 2), (5 и 6) и (7 и 8), равно нулю, то есть: $d_{1,2} = M_2 - M_2 = 0$; $d_{5,6} = M_6 - M_6 = 0$; $d_{7,8} = M_9 - M_9 = 0$. Поэтому если рассмотрим простое число с порядковым номером 1, то его ближайшим соседним простым числом с правой стороны является простое число с порядковым номером 2, а не число с порядковым номером 3 *(рис.2,с)*. Если рассмотрим простое число с порядковым номером 4, его ближайшим соседом с правой стороны является простое число с порядковым номером 5, а не число с порядковым номером 6. В этом легко убедиться, если разность между значениями этих чисел определим по формуле (1). Аналогичным образом легко можно определить, что ближайшим соседом простого числа с порядковым номером 9 с левой стороны является простое число с порядковым номером 8, а не число с порядковым номером 7, хотя эти два простого числа (7 и 9) расположены в одной строке.

Заметим, что числа с порядковыми номерами (1 и 2), (5 и 6) и (7 и 8) являются близнецами *(рис.2,с)*.

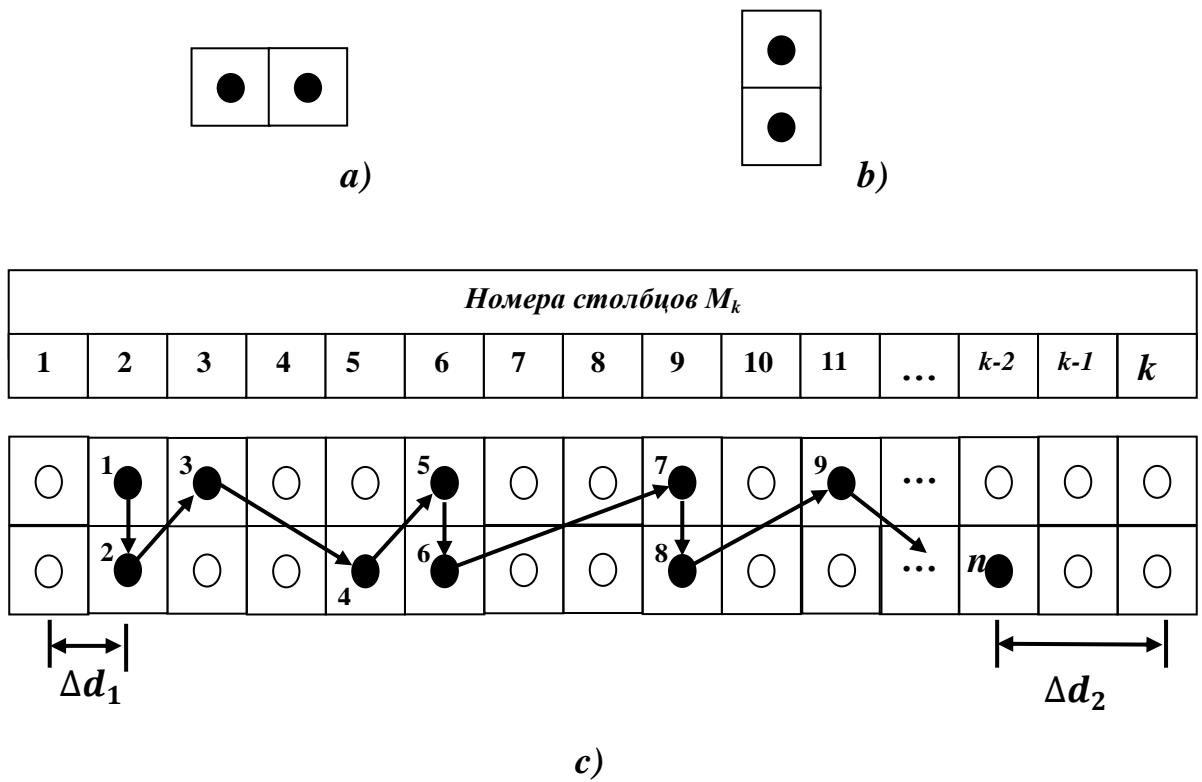

*Рис.2.* **Расположение простых чисел в паре строк-близнецов**

Здесь следует сказать, что первое простое число может занимать ячейку, которая будет находиться не в первом столбце рассматриваемого фрагмента *(рис.2,с)*. Аналогичным образом можно сказать, что ячейка, в которой располагается последнее простое число с порядковым номером $n$, также может находиться не в последнем столбце. Эти случаи в (9) не учтены, т.е. в выражении (9) отсутствуют параметры $\Delta d_1$ и $\Delta d_2$, где $\Delta d_1$ – количество столбцов от начала рассматриваемого фрагмента до ячейки первого простого числа,  $\Delta d_2$– количество столбцов от ячейки последнего простого числа до конца  рассматриваемого фрагмента *(рис.2,с)*. Чтобы учесть это рассмотрим сумму:

$$\Delta d = \Delta d_1 + \Delta d_2.$$

Здесь следует сказать, что значение рассматриваемого параметра $\Delta d$ сопоставимо и скорее всего, равно расстоянию по столбцам от ячейки последнего простого числа с порядковым номером $n$ до ячейки, где будет расположено следующее простое число с порядковым номером $n + 1$. Это число, как следует из рисунка 2-с, находится за пределами анализируемого фрагмента, но оно для -го простого числа будет ближайшим соседним простым числом с правой стороны. Из этого следует, что:

$$\Delta d \simeq d_{n,n+1}$$

Если при этом параметр $\Delta d$ добавим в сумму, которая находится в числителе выражения (9), то в точности получим «длину» рассматриваемого фрагмента *(рис.2,с)*:

$$\sum_{i=1}^{n-1} d_{i,i+1} + \Delta d = M_k$$

При этом количество слагаемых в числителе (9) будет на 1 больше, т.е. количество рассматриваемых простых чисел как бы растет на 1 и станет равным $n + 1$. С учетом этого выражение (9) приобретает следующий вид:

$$d_{\text{ср}} = \frac{\sum_{i=1}^{n-1} d_{i,i+1} + \Delta d}{(n+1) - 1} = \frac{M_k}{n} \qquad (10)$$

Теперь рассмотрим реальный случай – случай фрагмента матрицы $A_{k-1}$ с размерами $(m_{k-1}, M_{k-1})$. Здесь $m_{k-1}$ - количество пар строк-близнецов, $M_{k-1}$ - количество столбцов в данном фрагменте матрицы $A_{k-1}$, которая соответствует простому числу $p_{k-1}$. Здесь следует сказать, что значение $M_{k-1}$ должно быть достаточно большим, чтобы имело смысл заниматься статистикой. Например, должно быть:

$$M_k \gg p_k!' \qquad (11)$$

Пусть $\pi_{k-1}$ - общее количество всех простых чисел, находящихся во всех парах строк-близнецов рассматриваемого фрагмента матрицы $A_{k-1}$, $\pi_{k-1,l}$ - количество простых чисел, содержащихся в одной выбранной паре строк-близнецов с порядковым номером $l$. При этом общее количество пар строк-близнецов должно быть равно $m_{k-1}$, т.е. $l$ принимает значение от 1 до $m_{k-1}$, одним словом $l = 1 \div m_{k-1}$.

Тогда из (10) применительно для случая матрицы $A_{k-1}$ получим:
$$d_{k-1,l} = \frac{M_{k-1}}{\pi_{k-1,l}} \tag{12}$$
где $d_{k-1,l}$ - среднее расстояние <u>*по столбцам*</u> между ячейками, где располагаются соседние простые числа, находящиеся в одной паре строк-близнецов с порядковым номером $l$.

Заметим, здесь и в дальнейшем первый индекс рассматриваемого параметра (в данном случае это индекс $k-1$) будет соответствовать порядковому номеру рассматриваемой матрицы, а второй индекс этого параметра (в данном случае это индекс $l$) есть порядковый номер анализируемой пары строк-близнецов.

Выше нами было высказано предположение о том, что простые числа в парах строк-близнецов каждой новой матрицы должны располагаться более тесно, чем в парах строк-близнецов предыдущей матрицы. Для анализа и оценки этого предположения проанализируем значение $d_{k-1,l}$ по парам строк-близнецов рассматриваемого фрагмента.

Из работ Зигеля [6], [7], [8] и других ученых [9], [10], [11], [12] следует, что если постоянные (разность) разных арифметических прогрессии равны между собой, а также первый член и разность каждой арифметической прогрессии взаимно простые, то в этих прогрессиях простые числа распределяются одинаково и идентично. С другой стороны, как было показано выше, последовательность чисел, находящихся в любой строке любой пары строк-близнецов рассматриваемой матрицы $A_k$ образует арифметическую прогрессию с одинаковой разностью, равной $p_k!'$. Причем первый член и разность этих арифметических прогрессии взаимно простые. Эти прогрессии друг от друга отличаются только значением первого члена, а в остальном они идентичны. Это означает, что в пределах рассматриваемого фрагмента количество простых чисел в любой паре строк-близнецов примерно равно между собой, т.е.:
$$\pi_{k-1,\text{ср}} \simeq \pi_{k-1,1} \simeq \pi_{k-1,2} \simeq \cdots \simeq \pi_{k-1,l} \simeq \cdots \simeq \pi_{k-1,m_{k-1}}, \tag{13}$$
где $\pi_{k-1,\text{ср}}$ - среднее количество простых чисел, содержащихся в одной паре строк-близнецов рассматриваемого фрагмента матрицы $A_{k-1}$.

Тогда из этого и из (12) получим:
$$d_{k-1,\text{ср}} \simeq d_{k-1,1} \simeq d_{k-1,2} \simeq \cdots \simeq d_{k-1,l} \simeq \cdots \simeq d_{k-1,m_{k-1}},$$
или
$$d_{k-1,\text{ср}} = \frac{M_{k-1}}{\pi_{k-1,\text{ср}}} \tag{14}$$

где $d_{k-1,\text{ср}}$ - среднее расстояние *по столбцам* между ячейками, где располагаются соседние простые числа, находящиеся в одной паре строк-близнецов рассматриваемого фрагмента матрицы $A_{k-1}$.

С другой стороны, из (13) следует, что
$$\pi_{k-1,\text{ср}} = \frac{\sum_{j=1}^{m_{k-1}} \pi_{k-1,j}}{m_{k-1}} = \frac{\pi_{k-1}}{m_{k-1}} \qquad (15)$$

Здесь обратим внимание на следующее:
1. В отдельно выбранной строке любой пары строк-близнецов любой матрицы, за исключением матрицы $A_0$ и $A_1$ *(рис.1, $A_0$ и рис.1, $A_1$)*, как было показано выше, не может быть простых чисел-близнецов.
2. Если два простого числа, как было показано выше, располагаются в одном столбце в пределах одной пары строк-близнецов, то они являются близнецами. В этом случае расстояние *по столбцам* между ячейками, где располагаются эти два простого числа – близнецы, должно быть равно нулю.
3. Здесь рассматриваемое расстояние $d_{k,\text{ср}}$ *по столбцам* между ячейками, где находятся соседние простые числа, не следует путать с разностью, т.е. с расстоянием между значениями соседних простых чисел. Очевидно, что в одной ячейке рассматриваемого фрагмента любой матрицы может находиться только одно число. При этом разность значений соседних простых чисел, которые находятся в разных, но пусть примыкающих по строке друг к другу ячейках фрагмента, может принимать любое четное число.

Как было показано в случае доказательства теоремы 1 из (7) следует, что при преобразовании матрицы $A_{k-1}$ в матрицу $A_k$ каждая пара строк-близнецов матрицы $A_{k-1}$ порождает в новой матрице по $p_k - 2$ новых пар строк-близнецов и еще 2 одиночных строк. В результате каждая пара строк-близнецов матрицы $A_{k-1}$ создает $p_k - 1$ новых пар строк. Итого $m_{k-1}$ пар строк-близнецов порождают $m_k = m_{k-1}(p_k - 2)$ новых пар строк-близнецов и еще $m'_k = m_{k-1}$ пар из одиночных строк в новой матрице $A_k$.

Одним словом, множество простых чисел в количестве $\pi_{k-1} = m_{k-1} M_{k-1}/d_{k-1,\text{ср}}$, находящихся во всех парах строк-близнецов рассматриваемого фрагмента матрицы $A_{k-1}$, перераспределяется по $m_k + m'_k = m_{k-1}(p_k - 1)$ новым парам строк уже нового фрагмента, соответствующего новой матрице $A_k$. Тогда количество простых чисел, находящихся во всех вновь образованных $m_k$ парах строк-близнецов фрагмента матрицы $A_k$ определяется выражением:
$$\pi_k = \frac{m_k\, \pi_{k-1}}{m_k + m'_k} = \frac{m_k \pi_{k-1}}{m_{k-1}(p_k - 1)}$$

Следовательно, количество простых чисел, находящихся в одной выбранной паре строк-близнецов рассматриваемого фрагмента новой матрицы $A_k$, в среднем равно:
$$\pi_{k,cp} = \frac{\pi_k}{m_k} = \frac{\pi_{k-1}}{m_{k-1}(p_k-1)} \qquad (16)$$

С другой стороны при преобразовании фрагмента матрицы $A_{k-1}$ во фрагмент матрицы $A_k$ количество столбцов в новом фрагменте, как было показано выше, уменьшается в $p_k$ раз, т.е. $M_k = M_{k-1}/p_k$.

Тогда из (14), (15) и (16) получим среднее расстояние по столбцам $d_{k,cp}$ между ячейками, где располагаются соседние простые числа, находящиеся в одной паре строк-близнецов рассматриваемого фрагмента новой матрицы $A_k$:
$$d_{k,\text{ср}} = \frac{M_k}{\pi_{k,cp}} = \frac{m_{k-1}M_{k-1}}{\pi_{k-1}} * \frac{p_k - 1}{p_k} = d_{k-1,\text{ср}}\frac{p_k - 1}{p_k} \qquad (17)$$

Из этого следует, что с ростом порядкового номера $k$ матрицы $A_k$ среднее расстояние $d_{k,\text{ср}}$ между ячейками соседних простых чисел в любой паре строк-близнецов матрицы $A_k$ уменьшается непрерывно.

Поэтому из-за растущей тесноты из бесконечного количества простых чисел, а также идентичности их распределения в любой паре строк-близнецов матрицы $A_k$ простые числа-близнецы скорее «образуются», нежели в случае предыдущей матрицы $A_{k-1}$.

В частности, по причине того, что, как было показано выше, всеобщего перекоса простых чисел в парах строк-близнецов матрицы $A_3$ не существует, то его не будет и в парах строк-близнецов последующих матриц $A_4$, $A_5$, $A_6$, …, $A_k$, ... Стало быть, в каждой паре строк-близнецов любой матрицы $A_k$ имеются простые числа-близнецы.

**Теорема 2 доказана.**

Теперь рассмотрим проблему, поставленную перед данной работой.

**Теорема 3.** *Количество простых чисел-близнецов бесконечно.*

**Доказательство теоремы 3.**

Как было показано выше, каждая матрица $A_k$ соответствует определенному простому числу $p_k$. Известно, что количество простых чисел бесконечно, следовательно, количество разновидностей матриц $A_k$ тоже бесконечно. С другой стороны простые числа-близнецы могут находиться только в парах строк-близнецов.

Из теоремы 1 и из (7) следует, что по мере роста порядкового номера матрицы $A_k$ количество пар строк-близнецов в этой матрице непрерывно растет, причем
$$\lim_{k \to \infty} m_k = \lim_{k \to \infty} (p_k - 2)!' = \infty$$

Из теоремы 2 также следует, что в любой паре строк-близнецов имеются простые числа-близнецы. Это все означает, что количество простых чисел-близнецов бесконечно.

Этот вывод также неизбежно следует из (17).

Если в (17) параметр $d_{k-1,\text{ср}}$ выразим через $d_{k-2,\text{ср}}$ - среднее расстояние по столбцам между ячейками, где располагаются соседние простые числа в паре строк-близнецов матрицы $A_{k-2}$, то получим

$$d_{k-1,\text{ср}} = d_{k-2,\text{ср}} \frac{p_{k-1}-1}{p_{k-1}}$$

Тогда

$$d_{k,\text{ср}} = d_{k-2,\text{ср}} \frac{p_{k-1}-1}{p_{k-1}} * \frac{p_k-1}{p_k}$$

Далее аналогичным образом параметр $d_{k-2,\text{ср}}$ преобразуем через $d_{k-3,\text{ср}}$, затем $d_{k-3,\text{ср}}$ через $d_{k-4,\text{ср}}$, и т.д. преобразование продолжим до $d_{1,\text{ср}}$, тогда получим, что среднее расстояние по столбцам между ячейками простых чисел в любой паре строк-близнецов матрицы $A_k$ в общем виде определяется выражением:

$$d_{k,\text{ср}} = d_{1,\text{ср}} \frac{(p_1-1)(p_2-1)\ldots(p_k-1)}{p_1 p_2 * \ldots * p_k} = d_{1,\text{ср}} \frac{(p_k-1)!'}{p_k!'} \qquad (18)$$

где $k \geq 2$ и $d_{1,\text{ср}}$ – среднее расстояние между ячейками соседних простых чисел, находящихся в единственной неокрашенной строке матрицы $A_1$ *(рис.1, $A_1$-в)*.

Из (18) следует, что с ростом порядкового номера $k$ матрицы $A_k$ среднее расстояние $d_{k,\text{ср}}$ между ячейками простых чисел в любой паре строк-близнецов матрицы $A_k$ уменьшается непрерывно. Причем при $k \to \infty$ получим, что $d_{k,\text{ср}} \to 0$. В самом деле, проделав небольшое преобразование, из (18) получим, что

$$\ln d_{k,\text{ср}} = \ln d_{1,\text{ср}} - \left( \sum_{i=1}^{k} \frac{1}{p_i} + \sum_{i=1}^{k} \sum_{m=2}^{\infty} \frac{1}{m * p_i^m} \right)$$

Здесь бесконечный ряд, состоящий из обратных величин простых чисел, как было показано Эйлером, расходится [13], т.е.:

$$\sum_{i=1}^{\infty} \frac{1}{p_i} = \infty$$

Следовательно,

$$\lim_{k \to \infty} d_{k,\text{ср}} = 0 \qquad (19)$$

Это означает, что по мере роста порядкового номера матрицы $A_k$ среднее расстояние $d_{k,\text{ср}}$ между ячейками, где находятся соседние простые числа, стремиться к нулю. При этом бесконечность близнецов не только неизбежно, но и очевидно, так как из (7) следует, что при $k \to \infty$ количество пар строк-близнецов стремится к бесконечности: $\lim_{k \to \infty} m_k = \infty$, с другой стороны как следует из теоремы 2, в каждой из этих пар строк-близнецов имеются простые числа-близнецы.

**Теорема 3 доказана.**

***Резюме.***

*В работе предлагаются матрицы простых чисел, при помощи которых легко можно генерировать последовательность простых чисел. В статье также предлагается ряд теорем, при помощи которых доказывается бесконечность количества простых чисел-близнецов.*

***Ключевые слова****: простые числа, простые числа-близнецы, составные числа, натуральные числа, алгоритмы, арифметическая прогрессия, матрица простых чисел, специальный факториал, генерирование простых чисел.*